\documentclass[fleqn,12pt]{amsart}
\usepackage{enumerate,fullpage,nicefrac,amsmath,amssymb,multicol,amsthm, mathrsfs}
\usepackage{tikz}
\usetikzlibrary{arrows,shapes}

\newcommand{\comment}[1]{}

\newcommand{\conv}[1]{\text{conv}(#1)} 
\newcommand{\rbd}[1]{\text{rbd}(#1)}
\newcommand{\rint}[1]{\text{rint}(#1)}
\newcommand{\Ell}[0]{\mathscr{L}} 
\newcommand{\cdim}[1]{\text{cdim}(#1)}

\newtheorem{theorem}{Theorem}[section]
\newtheorem{corollary}[theorem]{Corollary}
\newtheorem{lemma}[theorem]{Lemma}

\newtheorem{proposition}[theorem]{Proposition}

\newtheorem{observation}[theorem]{Observation}

\theoremstyle{definition}

\theoremstyle{definition}

\allowdisplaybreaks

\begin{document}

\title[Order Dim. of Convex Geom.]{On the Order Dimension of Convex Geometries}

\author[J. Beagley]{Jonathan E. Beagley}

\address{Department of Mathematical Sciences\\ George Mason University, 4400 University Drive, Fairfax, VA 22030}

\email{jbeagley@gmu.edu}


\begin{abstract}
	We study the order dimension of the lattice of closed sets for a convex geometry.  Further, we prove the existence of large convex geometries realized by planar point sets that have very low order dimension.  We show that the planar point set of Erd\H{o}s and Szekeres from 1961 which is a set of $2^{n-2}$ points and contains no convex $n$-gon has order dimension $n-1$ and any larger set of points has order dimension strictly larger than $n-1$.
\end{abstract}

\maketitle

\section{Introduction}
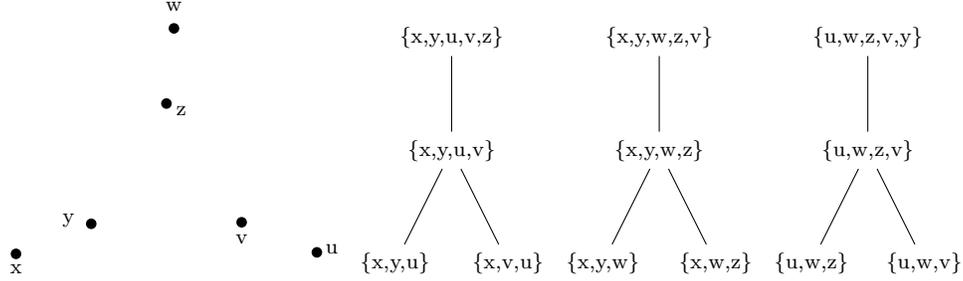
\begin{figure}
\begin{tikzpicture}
	\fill (0cm,0cm) circle (2pt);
	\fill (1cm,.4cm) circle (2pt);
	\fill (2cm,2cm) circle (2pt);
	\fill (2.1 cm,3cm) circle (2pt);
	\fill (3cm,.42cm) circle (2pt);
	\fill (4cm,.02cm) circle (2pt);
	
	\draw (0cm,-.2cm) node {x};
	\draw (0.7cm,.45cm) node {y};
	\draw (2.2cm ,1.9cm) node {z};
	\draw (2.1cm,3.3cm) node {w};
	\draw (3cm,.20cm) node {v};
	\draw (4.2cm,.06cm) node {u};
\end{tikzpicture}
\begin{tikzpicture}
\node {\{x,y,u,v,z\}}
child {node {\{x,y,u,v\}}
child {node {\{x,y,u\}}}
child {node {\{x,v,u\}}}
};
\end{tikzpicture}
\begin{tikzpicture}
\node {\{x,y,w,z,v\}}
child {node {\{x,y,w,z\}}
child {node {\{x,y,w\}}}
child {node {\{x,w,z\}}}
};
\end{tikzpicture}
\begin{tikzpicture}
\node {\{u,w,z,v,y\}}
child {node {\{u,w,z,v\}}
child {node {\{u,w,z\}}}
child {node {\{u,w,v\}}}
};
\end{tikzpicture}
\caption{A six point set and its poset of copoints}
\label{sixPointSet}
\end{figure}

Let $X$ be a finite set.  Let $\mathcal{C}$ be a collection of subsets of $X$ with the properties that; $\emptyset \in \mathcal{C}$, $X \in \mathcal{C}$, and for all $A,B \in \mathcal{C}$ then $A \cap B \in \mathcal{C}$.  $\mathcal{C}$ is called an \textit{alignment} of $X$ and when partially ordered by inclusion forms a lattice, $\mathcal{L} = (\mathcal{C},\subseteq)$.  For all $C \subseteq X$ we define $\Ell(C)$ to be the intersection of all $A \in \mathcal{C}$ such that $C \subseteq A$.  A subset $C$ of $X$ is \textit{closed} if $C = \Ell(C)$.  We say that $\Ell$ is \textit{anti-exchange} if given any $C \in \mathcal{C}$ and two distinct elements $p$ and $q$ in $X$, neither in $C$, then $q \in \Ell(C \cup p)$ implies that $p \notin \Ell(C \cup q)$.  A pair $(X, \Ell)$ for which $\Ell$ has the anti-exchange property is called a \textit{convex geometry} or \textit{antimatroid}.  Edelman and Jamison \cite{EJ1985} presented several equivalent definitions of convex geometry.  A closed subset of a convex geometry, $A \in \mathcal{C}$, is a \textit{copoint} if there is exactly one $B \in \mathcal{C}$ such that $|B \backslash A| = 1$.  Copoints are the meet-irreducible elements of the lattice $\mathcal{L}$.  The unique element in $B \backslash A$ is denoted $\alpha(A)$.  We say that the copoint $A$ is \textit{attached} to $\alpha(A)$.  The set of copoints partially ordered by inclusion is denoted $M(X)$.  We say $B \subseteq X$ is \textit{independent} if for all $p \in B$, $p \notin \Ell(B \backslash p)$. \\
\indent
The order dimension of a poset, $P = (X,\leq)$, sometimes called the dimension of $P$, is the least positive integer $t$ for which there exists a family $\mathcal{R} = \{L_1, L_2, \ldots, L_t\}$ of linear extensions of $P$ so that $P = \cap \mathcal{R}$.  Any family of linear extensions $\mathcal{R}$ such that $\cap \mathcal{R} = P$ is called a \textit{realizer} of $P$.  We specify order dimension to distinguish it from the convex dimension studied by Edelman and Saks \cite{ES1988}.  The convex dimension of $(X,\Ell)$ is the smallest number of chains needed to cover $M(X)$.  For the lattice of closed sets, $\mathcal{L}$, of convex geometry $(X,\Ell)$, we denote the order dimension by $\dim(\mathcal{L})$, convex dimension $\cdim{\mathcal{L}}$, and the size of the largest independent set $b(\mathcal{L})$. For any convex geometry, $\cdim{\mathcal{L}} \geq \dim(\mathcal{L}) \geq b(\mathcal{L})$ with strict inequalities possible (\cite{ES1988}).\\
\indent
As a natural example for a convex geometry, let $X$ be a set of $n$ points in $\mathbb{R}^2$ and for $A \subseteq X$ let $\Ell(A) = \conv{A} \cap X$.  A convex geometry is said to be \textit{realizable} in $\mathbb{R}^2$ if its lattice of closed sets is equivalent to the lattice of closed sets for some planar point set.  Figure \ref{sixPointSet} is a planar point set with its poset of copoints.  Morris \cite{morris2006} describes an algorithm to compute all the copoints of a planar point set in general position (no three on a line).  Start with a directed vertical line through $p \in X$.  Call the part of the line above $p$ the \textit{head} and the part of the line below $p$ the \textit{tail}.  Rotate the line clockwise around $p$, noting the order in which the points of $X \backslash p$ are met by the line.  If a point $q$ is met by the head of the line, write $q$, and if $q$ is met by the tail of the line, write $-q$.  The sequence of $2n-2$ symbols written as the line makes a complete revolution around $p$, viewed as a circular sequence, is called the \textit{circular local sequence} of $p$.  At one or more places in the circular local sequence of $p$ there will be an element $q$ followed by an element $-r$.  At such a place we can find a copoint.  Let $m$ be a line through $p$ of which $q$ and $r$ are on the same side.  Let $H$ be the open halfspace defined by $m$ that contains $q$ and $r$.  Then $H \cap X$ is a copoint attached to $p$. \\
\indent
Of particular interest to the author is a problem of Erd\H{o}s and Szekeres (\cite{erdosS1935},\cite{erdosS1961},\cite{MS2000}): For any $n \geq 3$, to determine the smallest positive integer $N(n)$ such that any set of at least $N(n)$ points in general position in the plane contains $n$ points that are the vertices of a convex $n$-gon.  It was established that $N(n) > 2^{n-2}$ and conjectured that $N(n) = 2^{n-2}+1$.  This conjecture has been verified for $n =$ 3,4,5, and 6.  We restate this conjecture in terms of convex geometries realized by points in $\mathbb{R}^2$, noting that an independent set of size $k$ in such a convex geometry corresponds to the vertex set of a convex $k$-gon.  For any convex geometry $(X,\Ell)$ realizable by a planar point set in general position with lattice of closed sets $\mathcal{L}$, it is conjectured that if $|X| > 2^{n-2}$  then $b(\mathcal{L}) \geq n$.  We show in section 4 of this paper that if $|X| > 2^{n-2}$ then $\dim(\mathcal{L}) \geq n$.
%
\tiny
\normalsize
%
\section{Critical Pairs and Copoints}
\begin{figure}
	\begin{tikzpicture}
	\node (1234) at (0,4) {1234};
	\node (123) at (-1,3) {123};
	\node (234) at (1,3) {234};
	\node (12) at (-2,2) {12};
	\node (23) at (0,2) {23};
	\node (24) at (2,2) {24};
	\node (1) at (-3,1) {1};
	\node (2) at (-1,1) {2};
	\node (3) at (1,1) {3};
	\node (0) at (0,0) {$\emptyset$};
	\draw (1234) -- (123)
				(1234) -- (234)
				(123) -- (12)
				(123) -- (23)
				(234) -- (23)
				(234) -- (24)
				(12) -- (1)
				(12) -- (2)
				(23) -- (2)
				(23) -- (3)
				(24) -- (2)
				(1) -- (0)
				(2) -- (0)
				(3) -- (0);
	\end{tikzpicture}
	\label{fourPoints}
	\caption{The lattice of closed sets for a non-atomic convex geometry}
\end{figure}
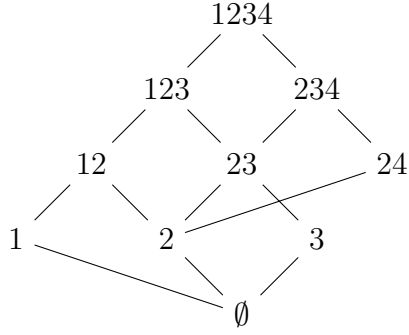

A \textit{critical pair} of elements in a poset, $P = (X, \leq)$ is an ordered pair $(A,B)$ where $A,B \in X$, $A$ is incomparable with $B$ and for all $U,D \in X$, $U > B$ implies that $U > A$ and $D < A$ implies that $D < B$.  Critical pairs are important in determining the order dimension of $P$.  It was shown by Rabinovitch and Rival \cite{RR1979} that $\mathcal{R}$ is a realizer of $P$ if and only if for every critical pair $(A,B)$, there is some $L \in \mathcal{R}$ for which $B < A$ in $L$.  We are able to conclude the following connection between copoints of a convex geometry and critical pairs of its lattice of closed sets.
\begin{theorem}\label{copointThm}
Let $(X,\Ell)$ be a convex geometry with lattice of closed sets $\mathcal{L} = (\mathcal{C},\subseteq)$.  For $A,B \in \mathcal{C}$, $(A,B)$ is a critical pair of $\mathcal{L}$ if and only if $B$ is a copoint attached to $p \in X$, $A = \Ell(p)$ and $A$ is incomparable with $B$.
\end{theorem}
\begin{proof}
	First, let $(A,B)$ be a critical pair of the lattice $\mathcal{L}$.  Suppose that $B$ is not a copoint, then there are distinct $p,q \in X \backslash B$ such that $B \cup p$ and $B \cup q$ are in $\mathcal{C}$.  Since $B \cup p$ and $B \cup q$ both properly contain $B$, they must also properly contain $A$.  The assumption that $p$ and $q$ are distinct implies that $B$ contains $A$, which contradicts $(A,B)$ being a critical pair.  Thus, $B$ must be a copoint attached to $p$.  Since every closed set $U$ that contains $B$ also contains $A$, we have that $B \cup p$ contains $A$.  Since $A$ is incomparable with $B$, $p \in A$.  Also, $\Ell(p) = A$ as any larger set containing $p$ properly contains $\Ell(p)$, and $B$ cannot contain $\Ell(p)$.\\
\indent
Suppose $B$ is a copoint attached to $p \in X$, $A = \Ell(p)$, and $A$ is incomparable with $B$.  Any $q \in A$, $q \neq p$ must be in $B$ because $B \cup p \in \mathcal{C}$ and $p \in B \cup p$.  Therefore, $D \subset A$ means that $D \subset B$.  Similarly, $U \supset B$ means that $U \supset A$.  Hence, $(A,B)$ is a critical pair of $\mathcal{L}$.
\end{proof}
In order to ensure that $(\Ell(p),B)$ is a critical pair, it is necessary for $\Ell(p)$ be incomparable with $B$ to eliminate the possibility where $\Ell(p) = B \cup p$.  The convex geometry with its lattice of closed sets being a chain is an example.  We call a convex geometry $(X,\Ell)$ \textit{2-edge-connected} if the Hasse diagram of the lattice of closed sets viewed as a graph has no cut edge.
\begin{proposition}
Let $(X,\Ell)$ be a convex geometry with lattice of closed sets $\mathcal{L}=(\mathcal{C},\subseteq)$.  Then every copoint $B$ of $(X,\Ell)$ is incomparable with $\Ell(\alpha(B))$ if and only if $\mathcal{L}$ is 2-edge-connected.
\end{proposition}
\begin{proof}
Let $B$ be a copoint of $(X, \Ell)$.  There is exactly one closed set $A$ greater than $B$ with $|A| = |B| + 1$, this is $A = B \cup \alpha(B)$.  $A = \Ell(\alpha(B))$ if and only if there is no smaller closed set containing $\alpha(B)$.  So, any closed set containing $\alpha(B)$ must be strictly greater than $B$ and because $B$ is a copoint all chains from $\emptyset$ to $X$ in $\mathcal{L}$ must contain $B$ and $A$.  The removal of this edge from the Hasse diagram of the lattice of closed sets from $B$ to $A$ will disconnect the Hasse diagram, so it is a cut edge.  Thus, the edge from $B$ to $A$ is not a cut edge if and only if $B$ is incomparable with $\Ell(\alpha(B))$.
\end{proof}
We can guarantee that $(X,\Ell)$ is 2-edge-connected when $|X| > 1$ and $(X,\Ell)$ is \textit{atomic}, that is for all $p \in X$, $\Ell(p) = p$.  This yields the following corollary.
\begin{corollary}
Let $(X,\Ell)$ be an atomic convex geometry of size greater than 1 with lattice of closed sets $\mathcal{L} = (\mathcal{C},\subseteq)$.  Let $A,B \in \mathcal{C}$, $(A,B)$ is a critical pair of $\mathcal{L}$ if and only if $B$ is a copoint attached to $A$.  That is, $B$ is a copoint and $A = \alpha(B)$.
\end{corollary}
\begin{proof}
The forward direction follows from Theorem \ref{copointThm} and the fact that $B$ does not contain $\alpha(B)$ and $\Ell(\alpha(B)) = \alpha(B)$, so $B$ and $\alpha(B)$ are incomparable.  The opposite direction is as in Theorem \ref{copointThm}.
\end{proof}
A convex geometry need not be atomic to be 2-edge-connected as we see in Figure \ref{fourPoints}.
%
%
\section{Order Dimension of Convex Geometries}
We consider the critical digraph, $\mathcal{D}(P)$ for the poset $P = (X, \leq)$ as defined by Trotter \cite{posetbook}.  This digraph has vertex set equal to the set of critical pairs of $P$ and there is a directed edge $(A,B) \rightarrow (C,D)$ anytime $C \leq B$.  In the case where $P = (\mathcal{C},\subseteq)$, the lattice of closed sets for convex geometry $(X,\Ell)$, this is simply where $C \subseteq B$ (and if $(X,\Ell)$ is atomic, $\alpha(D) \in B$).  The minimal cycles of this digraph induce a hypergraph, $\mathcal{H}^{C}_{P}$ with vertices from the set of critical pairs and hyperedges being minimal cycles of $\mathcal{D}(P)$.  For any cycle of $\mathcal{D}(P)$ there is no linear extension of $P$ that reverses all critical pairs belonging to that cycle.  Felsner and Trotter \cite{FT2000} showed that the chromatic number of this hypergraph is the order dimension of $P$ and greater than or equal to the chromatic number of the graph, $G^{C}_{P}$, induced from $\mathcal{H}^{C}_{P}$ by only considering the edges of size 2.  We state this as a lemma.
\begin{lemma}[\cite{FT2000}, 3.3]\label{FTLemma}
For every poset $P$, $\dim(P) = \chi(\mathcal{H}^{C}_{P}) \geq \chi(G^{C}_{P})$
\end{lemma}
For $(X,\Ell)$ a convex geometry with lattice of closed sets $\mathcal{L} = (\mathcal{C},\subseteq)$ we define $\mathcal{H}((X,\Ell)) = \mathcal{H}^{C}_{\mathcal{L}}$ as described above. We also define $\mathcal{G}((X,\Ell))$ with vertex set $M(X)$ and $AB \in E(\mathcal{G}((X,\Ell)))$ if and only if $\alpha(A) \in B$ and $\alpha(B) \in A$.  $\mathcal{G}((X,\Ell))$ was described and studied by Morris in \cite{morris2006} for convex geometries realizable by planar point sets in general position.  It was proven that a $k$-clique of $\mathcal{G}((X,\Ell))$ corresponds to the vertex set of convex $k$-gon in $X$ and the vertex set of a convex $k$-gon corresponds to a $k$-clique of $\mathcal{G}((X,\Ell))$.  The proof of this result can be extended to all convex geometries, where the $k$-cliques of $\mathcal{G}((X,\Ell))$ correspond to independent sets of size $k$. The correspondence of $k$-cliques of $\mathcal{G}((X,\Ell))$ and independent sets of size $k$ is not bijective.  We make note of the relationship between the graphs $\mathcal{G}((X,\Ell))$ and $G^{C}_{\mathcal{L}}$ for 2-edge-connected convex geometries.
\begin{theorem} \label{graphiso}
If $(X,\Ell)$ is 2-edge-connected with lattice of closed sets $\mathcal{L} = (\mathcal{C},\subseteq)$, then the graph $\mathcal{G}((X, \Ell))$ is isomorphic to the graph $G^{C}_{\mathcal{L}}$.
\end{theorem}
\begin{proof}
The bijection between the vertex sets is given by Theorem \ref{copointThm}, where $A$ maps to $(\Ell(\alpha(A)),A)$.  The edge of $AB \in E(\mathcal{G}((X,\Ell)))$ means $\alpha(B) \in A$ and $\alpha(A) \in B$.  This occurs if and only if $(\Ell(\alpha(A)),A) \rightarrow (\Ell(\alpha(B)),B)$ and $(\Ell(\alpha(B)),B) \rightarrow (\Ell(\alpha(A)),A)$ are directed edges in $\mathcal{D}(P)$.
\end{proof}
We make use of Theorem \ref{graphiso} and Lemma \ref{FTLemma} to yield a lower bound for the order dimension of 2-edge-connected convex geometries.
\begin{corollary} \label{myodim}
For any 2-edge-connected convex geometry $(X, \Ell)$ with lattice of closed sets $\mathcal{L} = (\mathcal{C},\subseteq)$, $\dim(\mathcal{L}) = \chi(\mathcal{H}^{C}_{\mathcal{L}}) \geq \chi(G^{C}_{\mathcal{L}}) = \chi(\mathcal{G}((X, \Ell)))$.
\end{corollary}
There are examples of posets $P$ where $\chi(\mathcal{H}^{C}_{P}) > \chi(G^{C}_{P})$ (\cite{posetbook}, \cite{FT2000}).  The remainder of this paper uses Theorem \ref{graphiso} and Corollary \ref{myodim} to calculate the order dimension of 2-edge-connected convex geometries.  To do so, it is necessary to first find and understand hyperedges of size greater than size 2 in $\mathcal{H}((X,\Ell))$.  Consider the convex geometry realized by the point set in Figure \ref{sixPointSet} with its poset of copoints.  One can quickly verify that both $(z, uwv), (v, xyu), (y,xwz)$ and $(z,xyw), (y, xvu), (v, uwz)$ are minimal cycles of length 3 for this convex geometry in its critical digraph.  However, $\chi(\mathcal{G}((X,\Ell))) = \chi(\mathcal{H}((X, \Ell))) =$ 4, so the order dimension is 4.\\
\indent
We define a cycle of copoints of length $l$, to be a collection $\mathcal{A}$ of $l$ copoints $(A_1, \ldots, A_{l})$ such that $\alpha(A_1) \in A_2, \alpha(A_2) \in A_3, \ldots, \alpha(A_{l-1}) \in A_l, \alpha(A_l) \in A_1$ and define $\alpha(\mathcal{A}) = \{\alpha(A_i): i =1, 2, \ldots, l\}$.  The following observations come directly from the anti-exchange property of a convex geometry and the assumption that the cycle is minimal.
\begin{observation}
	If $\mathcal{A}$ is a minimal cycle of length $l > 2$, then:
	\begin{itemize}
		\item
			$\alpha(A_i) \neq \alpha(A_j)$ unless $i = j$
		\item
			$\alpha(A_i) \notin A_j$ unless $j = i+1$ or $j = l$ and $i = 1$
		\item
			$A_i \not\subseteq A_j$ for $i \neq j$
		\item
			$\alpha(A_i)$ cannot be an extreme point of the convex geometry
	\end{itemize} 
\end{observation}
%

%
Let $X$ be a finite set of points in $\mathbb{R}^n$ and for all $A \subseteq X$, we define $\Ell(A) = \conv{A} \cap X$.  $(X,\Ell)$ is an atomic convex geometry and hence is 2-edge-connected, we say such a convex geometry is \textit{realizable} in $\mathbb{R}^n$.  If $A$ is a copoint of $(X,\Ell)$ attached to $\alpha(A)$, then $\alpha(A) \notin \conv{A}$, so $\conv{A}$ and $\alpha(A)$ can be properly seperated by a hyperplane $H_{A}$ in $\mathbb{R}^n$ (Theorem 2.4.10, \cite{convexity}) with $\alpha(A) \in H_{A}$ and $\conv{A}$ a subset of an open halfspace defined by $H_{A}$.  All copoints of $(X,\Ell)$ can be represented by such $H_{A}$.  We use this description of copoints to study $\alpha(\mathcal{A})$ where $\mathcal{A}$ is minimal cycle of $(X,\Ell)$.
\begin{proposition}
Let $(X,\Ell)$ be a convex geometry realizable by a set of points in $\mathbb{R}^n$.  If $\mathcal{A}$ is a minimal cycle of copoints with length $l \geq 2$, that is $\mathcal{A} \in E(\mathcal{H}((X,\Ell)))$ with $|\mathcal{A}| = l$, then $\alpha(\mathcal{A})$ is the vertex set of a polytope with $l$ vertices.
\end{proposition}
\begin{proof}
The result is trivial if $l = 2$, so let $l > 2$.  Suppose $\alpha(\mathcal{A})$ is not the vertex set of a polytope with $l$ vertices, then there is at least one $\alpha(A_i)$ such that $\alpha(A_i) \in \Ell(\{\alpha(A_j) : j \neq i\}) = C$.  By assumption, $\alpha(A_i) \in A_{i+1}$ so there is some $\lambda > 0$ such that $\lambda \alpha(A_i) + (1-\lambda) \alpha(A_{i+1}) = x \in \rbd{\conv{C}}$, the relative boundary of $\conv{C}$.  So, $x \in F$ a minimal face of $\conv{C}$ containing $x$, with $W \subset C$, $F = \conv{W}$.  Then, either $x \in W$ or $x \in \rint{F}$, the relative interior of $F$.  It is clear that any open halfspace with $\alpha(A_{i+1})$ on its boundary that contains $\alpha(A_i)$ must also contain $x$.  In particular, there is a hyperplane $H_{A_{i+1}}$ with $\alpha(A_{i+1})$ on its boundary with an open halfspace containing $\alpha(A_i)$ and hence $x$.  Thus, if $x \in W$ then $x \in A_{i+1}$ and $\mathcal{A}$ was not a minimal cycle.\\
\indent
Otherwise, $x \in \rint{F}$ and we write $x = \displaystyle\sum_{j=1}^{k} \lambda_j w_j$ for $w_j \in W$, $\lambda_j > 0$, and $\displaystyle\sum_{j=1}^{k} \lambda_j = 1$.  Let $H = \{y \in \mathbb{R}^n : y \cdot a = \alpha_0 \}$ be a hyperplane of $\mathbb{R}^n$ containing $\alpha(A_{i+1})$ and $x$ (and thus $\alpha(A_i)$).  So, $x \cdot a = \alpha_0 = (\displaystyle\sum_{j=1}^{k} \lambda_j w_j) \cdot a = \displaystyle\sum_{j=1}^{k} \lambda_j (w_j \cdot a)$.  Since $\lambda_j > 0$, if there is some $w_{j_1}$ such that $w_{j_1} \cdot a > \alpha_0$, then there is some and $w_{j_2} \cdot a < \alpha_0$.  In particular, there is a hyperplane $H_{A_{i+1}} = \{y \in \mathbb{R}^n: y \cdot a_{i+1} = \alpha_{i+1} \}$ with $\alpha(A_i) \cdot a_{i+1} < \alpha_{i+1}$ and hence $x \cdot a_{i+1} = (\displaystyle\sum_{j=1}^{k} \lambda_j w_j) \cdot a_{i+1} = \displaystyle\sum_{j=1}^{k} \lambda_j (w_j \cdot a_{i+1})< \alpha_{i+1}$.  Thus, there is a $w_j \in W$ such that $w_j \cdot a_{i+1} < \alpha_{i+1}$.  Therefore, the open halfspace that contains $\alpha(A_i)$ also contains $\alpha(A_j^{'}) = w_j \in W \subset C$ with $j^{'}$ not equal to $i$ or $i+1$.  This implies that $\mathcal{A}$ was not a minimal cycle, which is a contradiction.
\end{proof}
\section{The Erd\H{o}s-Szekeres Point Sets}
Let $X$ be a finite set of points in the plane, $\mathbb{R}^2$, in general position and for all $A \subseteq X$, we say $\Ell(A) = \conv{A} \cap X$.  Throughout this section, we will abbreviate the convex geometry $(X,\Ell)$ to be $X$ with this implied closure operator $\Ell$.  We call $B \subseteq X$ a \textit{block} in the circular local sequence of a point $p \in X$ if all $b \in B$ are consecutive and met by the same part of the line $\ell$ rotated around $p$ as described in Section 1.\\
\indent
Given two planar point sets $L$ and $M$, we define a \textit{composition} of $L$ and $M$ to be a point set of $L$ together with a translation of $M$ in which 
\begin{enumerate}
	\item
		every point of $M$ has greater first coordinate than the first coordinates of points of $L$,
	\item
		the slope of any line connecting a point of $L$ to a point of $M$ is greater than the slope of any line connecting two points of $L$ or two points of $M$.
\end{enumerate}
We easily see using the second condition, that the circular local sequence of $l \in L$ contains $-M$ and $M$ as blocks and that of $m \in M$ contains $-L$ and $L$ as blocks.\\
\indent
Erd\H{o}s and Szekeres \cite{erdosS1961} describe a construction of large point sets without large subsets in convex position.  For each pair of positive integers $(i,j)$, we define the point set $ES(i,j)$ that has $\binom{i+j}{i}$ points.  For all positive integers, $k$, we first define $ES(0,k)$ and $ES(k,0)$ to be singletons.  For $i \geq 1, j \geq 1$ define $ES(i,j)$ to be a composition of $ES(i-1,j)$ and $ES(i,j-1)$.\\
\indent
The extended Erd\H{o}s-Szekeres point set $XES(k)$ is a composition of $ES(0,k), ES(1, k-1), \cdots, ES(k,0)$ where the compositions are performed in order from left to right.  The number of points in $XES(k)$ is $2^{k}$ and it has convex dimension $2k-2$ for each positive integer $k$ (\cite{morris2006}).  The size of the largest independent set is $XES(k)$ is $k+1$, we will show that this is also the order dimension.
\begin{lemma} \label{contain}
Let $P$ be a planar point set such that $P$ is a composition of $L$ and $M$.  If $A$ is a copoint of $P$ attached to $\alpha(A) \in L$ that contains $m \in M$, then $M \subseteq A$.  Also, if $B$ is a copoint of $P$ attached to $\alpha(B) \in M$ that contains $l \in L$, then $L \subseteq B$.
\end{lemma}
\begin{proof}
We show that if $A$ is a copoint of $P$ attached to $\alpha(A) \in L$ that contains $m \in M$, then $M \subseteq A$ as the proof of the other statement is similar. $A$ is the intersection of $X$ with an open halfplane through $\alpha(A)$ that contains both $q, r \in X$, where $q$ and $-r$ are consecutive symbols in the circular local sequence of $\alpha(A)$.  $M$ is a block in the circular local sequence of $\alpha(A)$ and there exists $m \in M \cap A$.  Thus, all of $M$ is contained in this open halfplane and $M \subseteq A$.
\end{proof}
Let $T(P)$ be a rooted full binary tree where each node, $N_L$, of the tree corresponds to a point set $L$ with descendants $N_{M_1}$ and $N_{M_2}$ where $L$ is a composition of $M_1$ and $M_2$.  The root of $T(P)$ corresponds to the point set $P$, which we say is a composition of the leaves of the tree.  We define the sequence $P_L$ to be the unique sequence of nodes in $T(P)$ starting with $N_L$ followed by its ancestor nodes and terminating at the root.  We say that the point set $M$ is an \textit{ancestor} of $L$, or dually $L$ is a \textit{descendant} of $M$, if the node $N_M$ appears in the sequence $P_L$.\\
\indent
Let $p \in P$ and $p \in M$, a descendant of $P$, with $M$ being a composition of $M_1$ and $M_2$.  Suppose that $p \in M_1$, we know that the circular local sequence of $p \in M$ contains $M_2$ as a block.  The requirement that the the slope of any line in connecting $p$ to a point outside of $M$ must have greater slope than any line connecting two points of $M$ means that $M_2$ is a block in the circular local sequence for $p$ in $P$.  Similarly, if $p \in M_2$, the circular local sequence of $p$ in $P$ contains $M_1$ as a block.
\begin{lemma} \label{arbitraryContain}
Let $P$ be a composition of point sets $L_1, L_2, \ldots, L_n$.  If $A$ is a copoint of $P$ attached to $\alpha(A) \in L_i$ that contains $p \in L_j$ for some $j \neq i$, then $L_j \subseteq A$.
\end{lemma}
\begin{proof}
Suppose that $i < j$. There is at least one node common to the sequences $P_{L_i}$ and $P_{L_j}$, namely $N_P$.  So there is a first node, $N_M$, common to both paths.  $M$ is the composition of $M_1$ and $M_2$ where $M_1$ is an ancestor of $L_i$ but not $L_j$ and $M_2$ is an ancestor of $L_j$ but not $L_i$.  $M_2$ is a block in the circular local sequence of $\alpha(A)$ in $M$, so by the discussion before this lemma $M_2$ is a block in the circular local sequence of $\alpha(A)$ in $P$.  Since $A$ is a copoint of $P$ attached to $\alpha(A) \in M_1$ that contains $p \in M_2$, by Lemma \ref{contain}, $M_2 \subseteq A$ and $L_j \subseteq A$.  The proof is similar if $j > i$.
\end{proof}
We have proven a result stronger than what Lemma \ref{arbitraryContain} claims, that for all $L_j$ which are descendants of $M_2$, $L_j \subseteq A$.  We are now able to show a result on an arbitrary composition of point sets and the hyperedges of their associated hypergraph, $\mathcal{H}(P)$.
\begin{theorem} \label{pseudotheorem}
Let $P$ be the composition of point sets $L_1, L_2, \ldots, L_n$.  If $\mathcal{H}(P)$ contains a hyperedge $\mathcal{E}$ such that $|\mathcal{E}| > 2$, then $\alpha(\mathcal{E}) \subseteq L_i$ for some $i = 1, 2, \ldots, n$.
\end{theorem}
\begin{proof}
Suppose not, then there is an $\{A_1, \ldots, A_{l} \} = \mathcal{E} \in E(\mathcal{H}(P))$ such that $|\mathcal{E}| > 2$ where $\alpha(\mathcal{E})$ contains elements from $\{L_i : i \in I \}$ for $I \subseteq \{1, 2, \ldots, n\}$.  First, suppose that $\alpha(\mathcal{E})$ contains elements only from $L_i$ and $L_j$ for some $i \neq j$.  Let $\alpha(A_1) \in L_i$ and $\alpha(A_2) \in L_j$.  Since $\alpha(A_1) \in A_2$, $L_i \subseteq A_2$ by Lemma \ref{arbitraryContain}.  There is some $k \geq 2$ such that $\alpha(A_k) \in L_j$ and $\alpha(A_{k+1}) \in L_i$ and $\alpha(A_k) \in A_{k+1}$, so $L_j \subseteq A_{k+1}$.  If $A_{k+1} = A_1$, then $\alpha(A_2) \in A_1$ and $|\mathcal{E}| = 2$.  For $A_{k+1} \neq A_{1}$, we know that $\alpha(A_{k+1}) \in A_2$, so $\mathcal{E}$ is not a minimal cycle. \\
\indent
Therefore suppose that any such hyperedge has $|I| > 2$.  Consider all paths $P_{L_i}$ in $T(P)$ for $i \in I$.  These paths have at least one node in common, $N_P$, so there is a first common node to all paths $P_{L_i}$, call this node $N_M$.  $M$ is the composition of $M_1$ and $M_2$ with $L_i$ a descendant of $M_1$ if $i \in I_1$ and $L_i$ a descendant of $M_2$ if $i \in I_2$.  The sets $I_1$ and $I_2$ must be disjoint and non-empty, because $N_M$ is the first node in common to all paths.  By the proof of Lemma \ref{arbitraryContain}, if any copoint $A$ is attached to $\alpha(A) \in M_1$ that contains $p \in M_2$, $M_2 \subseteq A$.  So, let $\alpha(A_1) \in M_1$ and $\alpha(A_2) \in M_2$.  Then, there is some $k \geq 2$ such that $\alpha(A_k) \in M_2$ and $\alpha(A_{k+1}) \in M_1$ so $M_2 \subseteq A_{k+1}$.  This means that $\alpha(A_{k+1}) \in A_2$ and $\alpha(A_2) \in A_{k+1}$.  Thus, $\alpha(\mathcal{E})$ is contained in either $M_1$ or $M_2$ and $\mathcal{E}$ was not minimal.  Hence $\mathcal{E}$ was not in $E(\mathcal{H}(P))$, and the only remaining possibility is that $\alpha(\mathcal{E}) \subseteq L_i$ for some $i = 1,2, \ldots, n$.
\end{proof}
We apply Theorem \ref{pseudotheorem} to the point set $XES(k)$ to compute $\dim(XES(k))$.
\begin{corollary}
	$\mathcal{H}(XES(k)) \cong \mathcal{G}(XES(k))$.  Further, $\dim(XES(k)) = k+1$ for positive integers $k$.
\end{corollary}
\begin{proof}
Since $XES(k)$ is the composition of $2^k$ singletons, Theorem \ref{pseudotheorem} shows that if there is a hyperedge of size greater than 2 in $\mathcal{H}(XES(k))$, every copoint of the hyperedge must be attached to the same point.  This violates our observation on minimal cycles.  Thus, $\mathcal{H}(XES(k)) \cong \mathcal{G}(XES(k))$.\\
\indent
It was shown that $\chi(\mathcal{G}(XES(k))) = k+1$ by Morris \cite{morris2006}.  Therefore, $\chi(\mathcal{H}(XES(k))) = k+1$ and the order dimension of $XES(k)$ is $k+1$.
\end{proof}
We may also state this result for a more general class of point sets.
\begin{corollary}
If the point set $P$ is the composition of singletons, then $\mathcal{H}(P) \cong \mathcal{G}(P)$.
\end{corollary}
The Erd\H{o}s-Szekeres conjecture is a statement on the size of the planar point set in general position for a given size of the largest independent set.  We use a theorem of Morris \cite{morris2006} to come to a conclusion on the maximum size of a planar point set in general position with fixed order dimension.  Morris studied pseudoline arrangements, however planar point sets in general position are equivalent to stretchable pseudoline arrangements.
\begin{theorem}[\cite{morris2006}, 4.5] \label{morrissize}
If $L$ is a pseudoline arrangement and $\chi(\mathcal{G}(L)) = k$, then $|L| \leq 2^{k-1}$.
\end{theorem}
As a direct corollary to this theorem, we are able to bound the size of a planar point set in general position by a function of its order dimension.  This implies that any planar point set in general position of size greater than $2^{n-2}$ must have order dimension at least $n$.
\begin{corollary}
	If $P$ is a planar point set in general position and $\dim(P) = k$, then $|P| \leq 2^{k-1}$.
\end{corollary}
%
%
We pose some questions based on these results.  Is there some convex geometry for which $\chi(\mathcal{H}((X,\Ell))) > \chi(\mathcal{G}((X,\Ell)))$?  That is, does there exist some convex geometry with lattice of closed sets $\mathcal{L} = (\mathcal{C},\subseteq)$ with $\dim(\mathcal{L}) > \chi(\mathcal{G}((X,\Ell)))$?  The author has yet to construct an example of this.\\
\indent
As we have seen, $k$-cliques of $\mathcal{G}((X,\Ell))$ correspond to independent sets of size $k$, so the size of the largest independent set of $(X, \Ell)$ is $\omega(\mathcal{G}((X,\Ell)))$, the clique number.  Is $\frac{\chi(\mathcal{G}((X,\Ell)))}{\omega(\mathcal{G}((X,\Ell)))} \leq c$ for some constant $c$?  If $X$ is a finite set of planar points in general position then an independent set of size $k$ in $(X, \Ell)$ is the vertex sets of convex $k$-gon.  It has been conjectured by Walter Morris that $c \leq 2$ for convex geometries generated by such planar point sets.
\section{Acknowledgements}
The author wishes to thank his advisor, Walter Morris, for many helpful discussions and comments.
%

%
\end{document}